\documentstyle[11pt,amstex]{article}                                

\def\AFOUR{%
\setlength{\textheight}{9.0in}%
\setlength{\textwidth}{5.75in}%
\setlength{\topmargin}{-0.375in}%
\hoffset=-.5in%
\renewcommand{\baselinestretch}{1.17}%
\setlength{\parskip}{6pt plus 2pt}%
}

\AFOUR                                           

\makeatletter
\def\section{\@startsection {section}{1}{\z@}{-3.5ex plus -1ex minus
 -.2ex}{2.3ex plus .2ex}{\large\bf}}
\def\subsection{\@startsection{subsection}{2}{\z@}{-3.25ex plus -1ex minus
 -.2ex}{1.5ex plus .2ex}{\normalsize\bf}}
\makeatother

\makeatletter
\@addtoreset{equation}{section}

\makeatother

\newcommand{\nc}{\newcommand}
\newcommand{\rnc}{\renewcommand}

\nc{\bea}{\begin{eqnarray}}
\nc{\eea}{\end{eqnarray}}
\nc{\be}{\bea}
\nc{\ee}{\eea}

\rnc{\a}{\alpha}
\nc{\ab}{\bar{\a}}
\nc{\ap}{\a^{+}}
\nc{\abm}{\ab^{-}}
\rnc{\b}{\beta}
\nc{\bb}{\bar{\b}}
\nc{\bbp}{\bb_{\zb}^{+}}
\nc{\bm}{\b_{z}^{-}}
\nc{\oa}{\overline{\a}}
\nc{\ob}{\overline{\b}}
\rnc{\gg}{\gamma}
\rnc{\d}{\delta}
\nc{\f}{\phi}
\nc{\fb}{\bar{\phi}}
\nc{\vf}{\varphi}
\nc{\p}{\psi}

\rnc{\c}{\chi}
\nc{\la}{\lambda}
\nc{\m}{\mu}
\nc{\n}{\nu}
\rnc{\o}{\omega}
\nc{\Om}{\Omega}
\rnc{\t}{\theta}
\nc{\eps}{\epsilon}
\rnc{\S}{\Sigma}
\nc{\F}{\Phi}

\nc{\trac}[2]{{\textstyle\frac{#1}{#2}}}

\nc{\ex}[1]{\mbox{e}^{\,\textstyle#1}}

\nc{\mat}[4]{\left(\begin{array}{cc}#1&#2\\#3&#4\end{array}\right)}

\nc{\som}[9]{\left(\begin{array}{ccc}#1&#2&#3\\#4&#5&#6\\#7&#8&#9%
\end{array}\right)}

\nc{\tr}{\mathop{\mbox{tr}}\nolimits}
\nc{\ad}{\mathop{\mbox{ad}}\nolimits}
\nc{\Tr}{\mathop{\mbox{Tr}}\nolimits}
\nc{\ST}{\mathop{\mbox{STr}}\nolimits}
\nc{\Det}{\mathop{\mbox{Det}}\nolimits}
\nc{\rk}{\mathop{\mbox{rk}}\nolimits}
\nc{\ra}{\rightarrow}
\nc{\Ra}{\Rightarrow}
\nc{\LRa}{\Leftrightarrow}
\nc{\ot}{\otimes}
\rnc{\ss}{\subset}
\nc{\nul}{\noindent\underline}
\nc{\non}{\nonumber\\}

\nc{\subs}[1]{{\vspace*{0.5cm}}%
{\noindent\underline{#1}}{\addcontentsline{toc}{subsection}{#1}}%
{\vspace*{0.3cm}}}

\nc{\zb}{\bar{z}}
\rnc{\lg}{\frak{g}}
\nc{\lt}{\frak{t}}
\nc{\lk}{\frak{k}}
\nc{\lh}{\frak{h}}
\nc{\pik}{\Pi_{\lk}}
\nc{\pip}{\Pi_{+}}
\nc{\pim}{\Pi_{-}}
\nc{\pih}{\Pi_{\lh}}
\nc{\jz}{J_{z}}
\nc{\jzh}{\jz^{\lh}}
\nc{\jzp}{\jz^{+}}
\nc{\jzm}{\jz^{-}}
\nc{\del}{\partial}
\nc{\dz}{\del_{z}}
\nc{\dzb}{\del_{\bar{z}}}
\nc{\az}{A_{z}}
\nc{\azb}{A_{\bar{z}}}
\nc{\g}{g^{-1}}
\nc{\dw}{\Delta_{W}}
\nc{\Ad}{{\mbox{Ad}}}
\nc{\ks}{Ka\-za\-ma-\-Su\-zu\-ki}
\nc{\KS}{\ks}
\nc{\ksm}{\ks\ model}
\rnc{\AA}{{\Bbb A}}
\nc{\BB}{{\Bbb B}}
\nc{\CC}{{\Bbb C}}
\nc{\PP}{{\Bbb P}}
\nc{\cpm}{\CC\PP(m)}
\nc{\cpn}{\CC\PP(n)}
\nc{\cp}[1]{\CC\PP(#1)}
\nc{\gmn}{G(m,m+n)}
\nc{\gmnk}{\gmn_{k}}
\nc{\cO}{{\cal O}}
\nc{\bcO}{\bar{\cO}}
\nc{\bO}{\bar{O}}
\nc{\oQ}{\overline{Q}}

\newtheorem{theorem}{Theorem}[section]

\newtheorem{corollary}[theorem]{Corollary}

\begin{document}
\global\parskip=4pt

\makeatletter
\begin{titlepage}
\begin{center}

\begin{flushright}
IC/99/58
\end{flushright}
\begin{center}
\vskip .5in

{\LARGE\bf A Geometric Interpretation of the $\chi_{y}$ Genus on
Hyper-K\"{a}hler Manifolds}\\

\vskip 0.4in
{\bf George Thompson}\footnote{email: thompson@ictp.trieste.it}
\vskip .1in
ICTP \\
P.O. Box 586 \\
34100 Trieste \\
Italy\\

\end{center}
\vskip .4in
\begin{abstract}
The group $SL(2)$ acts on the space of cohomology
groups of any hyper-K\"{a}hler manifold $X$. The $\chi_{y}$ genus of
a hyper-K\"{a}hler $X$ is shown to have a geometric interpretation as
the super trace of an element of $SL(2)$. As a by product one learns
that the generalized Casson invariant for a mapping torus is
essentially the $\chi_{y}$ genus.
\end{abstract}
\end{center}
\end{titlepage}
\makeatother
\begin{small}
\tableofcontents
\end{small}

\setcounter{footnote}{0}

\section{Introduction}
The $\chi_{y}$ genus of Hirzebruch is a very interesting and rather
powerful invariant. There are three significant values for $y$. At
$y=-1$ the $\chi_{y}$ genus is the Euler
characteristic, at $y=0$ it is the Todd genus while at $y=1$ it is the
signature. There seems to be, however, no geometric
understanding of the genus away from these prefered values of $y$. In
this short note, I prove that for (compact) hyper-K\"{a}hler manifolds, there
is, in fact, quite a clear geometric meaning to the genus. 

For
hyper-K\"{a}hler manifolds there is a natural $SL(2)$ action,
associated with the holomorphic 2-form, on the cohomology groups
$\bigoplus_{p} {\mathrm H^{q}}\left( X, \Omega_{X}^{p}\right)$ which
preserves $q$ and shifts $p$ by even integers. This means that
$(-1)^{q+p}$ is preserved. One can, therefore, take the graded trace of
an $SL(2)$ element, with the grading given by $(-1)^{p+q}$. Denote the
graded trace of $U \in SL(2)$ by $\ST U$.

The geometric meaning of the $\chi_{y}$ genus for hyper-K\"{a}hler $X$
is the content of the following
\begin{theorem}{Let $X$ be an irreducible compact hyper-K\"{a}hler
manifold of real
dimension $4n$. Let $U \in SL(2)$ and $y$ an eigenvalue of $U$, in
the two dimensional representation, then} \label{1}

\end{theorem}
\be
\ST \, U \, = \, \frac{\chi_{-y}}{y^{n}} . \label{id}
\ee
{\bf Remarks:}

1) Note that, since $h^{(p,q)}= h^{(2n-p,q)}$, the right hand side is
invariant under $y \rightarrow 1/y$ so that it does not depend on
which eigenvalue one picks. 

2) Once one expects that a result of this kind is true the proof turns
out to be embarrassingly easy.

The motivation for this result comes from the study of 3-manifold
invariants. Rozansky and Witten \cite{RW} indicated how, given a
hyper-K\"{a}hler manifold $X$, one could associate to the Mapping Torus
$T_{U}$, the invariant $\ST \, U \,$. In \cite{T}, I showed that one
could perform the associated path integral. The solution found
there is, in fact, the Riemann-Roch formula for the $\chi_{y}$ genus divided
by $y^{n}$. This motivated the above theorem, which can be proven
without recourse to physics. However, one can now read the derivation
in \cite{T} as a path integral proof of the Riemann-Roch formula for
the $\chi_{y}$ genus.

That path integral calculation of $\ST U$ gave, 
\be
\int_{X} \, {\mathrm
Todd}\left(TX_{{\Bbb C}}\right) \, \Det{\left(U\otimes I - I \otimes \ex{R}
\right) }^{1/2} . 
\ee
Which can be re-written as
\be
\int_{X} \, {\mathrm
Todd}\left(TX_{{\Bbb C}}\right) \, \prod_{i=1}^{n}\left(t-
2\cosh{x_{i}} \right) ,\label{mtf}
\ee
where $t$ is the character of $U$ in the 2-dimensional representation.
The $\chi_{y}$ genus is given by Riemann-Roch as \cite{NR}
\be
\chi_{-y} (X)= \int_{X} \, {\mathrm
Todd}\left(TX_{{\Bbb C}}\right) \prod_{i=1}^{2n}\left( 1- y
\ex{-x_{i}} \right) ,
\ee
but since $X$ is hyper-K\"{a}hler one has that $x_{i+n}=-x_{i}$ for
$i\leq n$. This means that
\be
\chi_{-y} (X)= \int_{X} \, {\mathrm
Todd}\left(TX_{{\Bbb C}}\right) \prod_{i=1}^{n}\left( (1+y^{2}) -
2y \, {\mathrm cosh} (x_{i}) \right) ,
\ee
so that this suggests (\ref{id}) on setting $t y = 1 +
y^{2}$. 

Consequently we have, in the notation of \cite{T}, 
\begin{corollary}{The Rozansky-Witten invariant $Z_{X}^{RW}[T_{U}] =
\chi_{-y}/y^{n}$, for $U\in SL(2, {\Bbb Z})$.}\end{corollary}

{\bf Further Remarks:}

1) The essential feature used here is the $SL(2)$ action that is made 
available by the holomorphic 2-form. Hence this is not the same as
thinking of $X$ as a K\"{a}hler manifold and making use of the usual
$SL(2)$ action that comes from the symplectic 2-form (Lefschetz
decomposition). 

2) There is a rather more general formula that was suggested by the
work of \cite{RW}. If one considers a ``mapping Riemann surface'', for
a Riemann surface, $\Sigma$, of genus $g$, then the Rozansky-Witten invariant
$Z^{RW}_{X}[\Sigma_{U}] = \ST U$ where $U\in {\mathrm Sp}(g)$ and
this group acts on $\bigoplus {\mathrm H}^{q}\left( X, \left
( \Omega_{X}^{*} \right)^{\otimes g} \right)$. In \cite{T} a
Riemann-Roch formula for this super trace was given which looks like a
Riemann-Roch formula for a generalized $\chi_{y}$ genus. That suggests
that the corresponding generalized $\chi_{y}$ can be rigorously shown
to be the super trace. This has important implications for 3-manifold
invariants. 

3) Similar, though not identical, path integral formulae are available
for general holomorphic symplectic manifolds.

4)Justin Sawon \cite{S} has made use of the weight system in \cite{RW} in an
ingenious way to get constraints on the Chern numbers of $X$.

\section{The Sl(2) Action on X}
The $SL(2, {\Bbb C})$ action on the cohomology groups of $X$, that we are
interested in, is perhaps best
explained at the level of the Lie algebra, ${\mathbf Lie}$
$SL(2) := sl(2)$. Let $L_{\eps}:
{\mathrm H^{q}}\left( X, \Omega_{X}^{p}\right) \rightarrow {\mathrm
H^{q}}\left( X, \Omega_{X}^{p+2}\right) $ be the map given by the
cup-product with the holomorphic 2-form $\eps$. Let $
\imath_{\eps}:{\mathrm  H^{q}}\left( X, \Omega_{X}^{p}\right)
\rightarrow  {\mathrm H^{q}}\left( X, \Omega_{X}^{p-2}\right)$
be contraction with respect to $\eps$. To fix conventions we note that
in local holomorphic coordinates if $\o \in \Omega^{(p,q)}(X)$, then,
suppressing the anti-holomorphic factors, (the Einstein summation
convention is in force)
\be
\o = \o_{I_{1}, \dots, I_{p}} \, dz^{I_{p}} \wedge \dots \wedge
dz^{I_{1}}, 
\ee
and
\be
\imath_{\eps}\o = \frac{p(p-1)}{2} \, \o_{I_{1}, I_{2}, I_{3}, \dots, I_{p}} \,
\eps^{I_{1}  I_{2}} \, dz^{I_{3}} \wedge \dots \wedge dz^{I_{p}} . 
\ee

The algebra satisfied by these operators is, by a straightforward
computation,
\be
\left[ \imath_{\eps} , L_{\eps}\, \right] = (n-p)
\ee
understood as a map ${\mathrm  H^{q}}\left( X, \Omega_{X}^{p}\right)
\rightarrow  {\mathrm H^{q}}\left( X, \Omega_{X}^{p}\right)$. The
generators of $sl(2)$ are then realized as
\be
\left(\begin{array}{cc}
0 & 1 \\
0 & 0 
\end{array}\right) \sim L_{\eps}  \; \; 
\left(\begin{array}{cc}
0 & 0 \\
1 & 0 
\end{array}\right) \sim \imath_{\eps}  \; \; 
\left(\begin{array}{cc}
1 & 0 \\
0 & -1 
\end{array}\right) \sim (n-p) .
\ee

The following is taken from the survey by Huybrechts \cite{H} (but see
also the original work by Fujiki \cite{F}). 
Let,
\be
{\mathrm H}^{q}\left(X, \Omega^{p}_{X}\right)_{\eps} := {\mathrm
ker} \, L^{n-p+1}_{\eps} , 
\ee
then the Lefschetz decomposition theorem tells
us that
\be
{\mathrm H}^{q}\left(X, \Omega^{p}_{X}\right) \bigoplus_{(p-l) \geq
{\mathrm max}(p-n,0)}  L^{p-l}_{\eps}{\mathrm H}^{q}\left(X,
\Omega^{2l-p}_{X}\right)_{\eps} .
\ee
One thinks of $L_{\eps}$ as a raising operator, and the ${\mathrm
H}^{q}\left(X, \Omega^{p}_{X}\right)_{\eps}$, for $0 \leq p\leq n$, are
the highest weight vectors of the $n-p+1$ dimensional irreducible
representations of $SL(2, {\Bbb C})$. One also has, by a straightforward count,
that
\be
{\mathrm dim}_{{\Bbb R}} {\mathrm H}^{q}\left(X,
\Omega^{p}_{X}\right)_{\eps} := h_{\eps}^{(p,q)} = h^{(p,q)}-h^{(p-2,q)}.
\ee

\section{Proof of Theorem \ref{1}}
The proof is by direct computation.

Let $t_{r}$ be the character
of $U$ in the
$r$ dimensional irreducible representation of $SL(2, {\Bbb C})$ and set
$t_{2}=t$. Note that $t_{1} = 1$, and I use the convention that
$t_{r}=0$ for $r\leq 0$, as well as $h^{(p,q)}=0$ if $p<0$. Then
\be
\ST U = \sum_{q=0}^{2n} \sum_{p=0}^{n} (-1)^{p+q} \, t_{n-p+1} \,
h_{\eps}^{(p,q)}.
\ee
One can re-write this expression as
\be
\ST U = \sum_{q=0}^{2n} \sum_{p=0}^{n} (-1)^{p+q} \, h^{(p,q)}
\left(t_{n-p+1}-t_{n-p-1} \right) . \label{stut}
\ee

Now notice that, on making use of Serre duality, which implies that
$h^{(p,q)}= h^{(2n-p,q)}$, that the $\chi_{y}$
genus satisfies,
\be
\frac{\chi_{-y}}{y^{n}} = \sum_{q=0}^{2n}\sum_{p=0}^{n-1} (-1)^{p+q}
h^{(p,q)} \left( y^{p-n} + y^{n-p} \right) \; + \sum_{q=0}^{2n}
(-1)^{q} h^{(n,q)} . \label{stuy}
\ee

A comparison of (\ref{stut}) and (\ref{stuy}) shows us that they agree
if we can set
\be
t_{r+1} - t_{r-1} = y^{r}+y^{-r} \;\;\; r > 0 . \label{eig}
\ee
For $r=1$ this reads as
\be
t y = y^{2} +1 ,
\ee
which is simply the characteristic polynomial for the two-dimensional
representation of $U$, where $y$ is an eigenvalue and $t$ is the
trace. We make this identification, then (\ref{eig}) is a standard
relationship between characters and eigenvalues for $SL(2)$.\begin{flushright}
$\Box$
\end{flushright}

\subsubsection*{Acknowledgments}
I would like to thank M. Blau, L. G\"{o}ttsche and I. King for
discussions. Special thanks are due to M. S. Narasimhan who made the right
observations and the right remarks at the right time.

\rnc{\Large}{\normalsize}

\end{document}